\input amstex.tex 
\document
\magnification=\magstep1
\baselineskip 18pt
\NoBlackBoxes

\par \vskip 2pc \noindent
$$\text{\rm {\bf NIEMEIER LATTICES AND K3 GROUPS}}$$

\par \vskip 1pc \noindent
$$\text{\rm Dedicated to Professor I. Dolgachev on the occassion of his sixtieth birthday}$$

\par \vskip 2pc \noindent
$$\text{\rm {\bf D. -Q. Zhang}}$$

\par \vskip 2pc \noindent
{\bf Abstract.}
In this note, we consider $K3$ surfaces $X$ with an action by the
alternating group $A_5$. We show that if a cyclic extension $A_5 . C_n$
acts on $X$ then $n = 1$, $2$, or $4$. We also determine the
$A_5$-invariant sublattice of the K3 lattice and its discriminant form.

\par \vskip 2pc \noindent
{\bf Introduction}

\par \vskip 1pc \noindent
We work over the complex numbers field ${\bold C}$.
A {\bf K3} surface $X$ is a simply connected projective surface
with a nowhere vanishing holomorphic 2-form $\omega_X$.
In this note, we will consider finite groups in $\text{Aut}(X)$.
An element $h \in \text{Aut}(X)$ is {\bf symplectic}
if $h$ acts trivially on the 2-form $\omega_X$.
A group $G_N \subseteq \text{Aut}(X)$ is {\bf symplectic}
if every element of $G_N$ is symplectic.

\par \vskip 1pc \noindent
According to Nikulin [Ni1], Mukai [Mu1] and Xiao [Xi], there are 
exactly 80 abstract finite groups which can act symplectically
on $K3$ surfaces.
Among these 80, there are exactly three non-abelain simple groups 
$A_5, \, L_2(7)$ and $A_6$.

\par \vskip 0.5pc \noindent
To be more precise, as in {\bf (1.0)} below, for every
finite group $G$ acting on a $K3$ surface $X$,
the symplectic elements of $G$ (i.e., those $h$
acting trivially on the non-zero 2-form $\omega_X$)
form a normal subgroup $G_N$ such that $G/G_N \cong \mu_I$
(the cyclic group of order $I$ in ${\bold C}^*$). 
Namely, we have $G = G_N . \mu_I$
(see {\bf Notation} below).
The natural number $I = I(G)$ is determined by $G$ and
called the {\bf transcendental value} of $G$.

\par \vskip 1pc \noindent
It is proved in [OZ3] and [KOZ1, 2] that when $G_N$ is either one
of the two bigger simple groups above,
the transcendental value $I(G) \ne 3$.
As expected or unexpected, the same is true for the smaller (indeed the smallest
non-abelian simple group) $A_5$:

\par \vskip 1pc \noindent
{\bf Theorem A.} Suppose that $G = A_5 . \mu_I$ acts faithfully on a $K3$ surface
(assuming that $G_N = A_5$).
Then $G = A_5 : \mu_I$ (semi-direct product) and $I = 1$, $2$, or $4$.

\par \vskip 1pc \noindent
In general, for a group of the form $G = A_5 . C_n$ acting on a $K3$ surface
(here $G_N$ might be bigger than $A_5$; and $C_n$ an abstract cyclic group of order $n$), 
we have a similar result:

\par \vskip 1pc \noindent
{\bf Theorem B.} Suppose that a group of the form $G = A_5 . C_n$ acts 
faithfully on a $K3$ surface. Then $G = A_5 : C_n$ and $I = 1$, $2$, or $4$.
Moreover, $G_N = A_5$ (and hence $C_n = \mu_n$ 
in the notation above or {\bf (1.0)}) unless $G = G_N = S_5$.
 
\par \vskip 1pc \noindent
We can determine the $A_5$-invariant sublattice of the $K3$ lattice 
in the result below, which has application in helping determine the
transcendental lattice $T_X$ and hence the surface itself (when 
rank $T_X = 2$).

\par \vskip 1pc \noindent
{\bf Theorem C.} Suppose that $A_5$ acts faithfully on a $K3$ surface $X$.
Then we have:

\par \noindent
(1) The $A_5$-invariant sublattice $L^{A_5}$ of the $K3$ lattice
$L = H^2(X, {\bold Z})$ has rank $4$.
The $A_5$-invariant sublattice $S_X^{A_5}$ of the Neron Severi lattice $S_X$
has rank equal to $1$ or $2$. 

\par \noindent
(2) The discriminant group 
$A_{L^{A_5}} = \text{\rm Hom}(L^{A_5}, {\bold Z})/L^{A_5}$ equals one of the following
(see Theorem {\bf (2.1)} for the corresponding intersection forms):
$$\gather
{\bold Z}/(30) \oplus {\bold Z}/(30), \hskip 0.5pc
{\bold Z}/(30) \oplus {\bold Z}/(10), \hskip 0.5pc
{\bold Z}/(20) \oplus {\bold Z}/(20), \\
{\bold Z}/(60) \oplus {\bold Z}/(20), \hskip 0.5pc
{\bold Z}/(60) \oplus {\bold Z}/(20) \oplus {\bold Z}/(2) \oplus {\bold Z}/(2).
\endgather$$

\par \vskip 1pc \noindent
{\bf Remark D.} 

\par \noindent
(1) In [Z2], it is proved that there is no faithful action of $A_5 . \mu_4$
on a $K3$ surface. So the $I$ in Theorems A and B can only be 1 or 2.

\par \noindent
(2) Theorem C is used in [Z2, Lemma 3.5]. The proofs of Theorems A, B and C
here are independent of the paper [Z2].

\par \noindent 
(3) Theorem C is applicable in the following situation:
Suppose in addition that a non-symplectic involution $g \in \text{\rm Aut}(X)$ commutes with
every element in $A_5$ and that the fixed locus $X^g$ 
is a union of a genus $\ge 2$ curve $C$ and $s$ ($\ge 1$) smooth rational curves $D_i$.
Then $S_X^{A_5}$ contains $L_0 = {\bold Z}[C, \sum_{i=1}^s D_i]$
as a sublattice of finite index $d_1$. Note that $L^{A_5}$ contains
$S_X^{A_5} \oplus T_X$ as a sublattice of finite index $d$.
So $|L_0| |T_X| = d_1^2 |S_X^{A_5}||T_X| = d_1^2 d^2 |L^{A_5}|$
while $-|L^{A_5}| = 30^2, 3 \times 10^2, 20^2$, $3 \times 20^2$, 
or $3 \times 40^2$ as given in Theorem C.
This is a restriction imposed on $|T_X|$. 
In [Z2], we determined $d_1$, $d$, and $|T_X|$ using the existence of 
the extra $\mu_4$ in (the impossible case:) $A_5 . \mu_4$
where $T_X$ then has the intersection form $\text{\rm diag}[2m, 2m]$ for some $m \ge 1$.

\par \noindent
(4) The same construction in [OZ3, Appendix] shows that there is a
smooth non-isotrivial family of $K3$ surfaces $f : {\Cal X} \rightarrow {\bold P}^1$
such that all fibres admit $A_6$ actions and infinitely many of them
are algebraic K3 surfaces. So, the symplectic part alone
can not determine the surface uniquely, and the study of transcendental value is needed.

\par \vskip 1pc \noindent
The main tools of the paper are the Lefschetz fixed point formula
(both the topological version and vector bundle version due to Atiyah-Segal-Singer [AS2, 3]),
the representation theory on the $K3$ lattice and the study on automorphism groups
of Niemeier lattices (in connection with Golay binary or ternary codes)
where the latter is much inspired by Conway-Sloane [CS], Kondo [Ko1] and 
Mukai [Mu2]. 

\par \vskip 1pc \noindent
We believe that the way of combining different very powerful machinaries to deduce results as done in the paper should be applicable
to the study of other problems. Our humble paper also demonstrates
the powerfulness and depth of these algebraic results in the study of geometry.

\par \vskip 1pc \noindent
{\bf Notation.}

\par \noindent
{\bf 1.} $S_n$ is the symmetric group in $n$ letters,
$A_n$ ($n \ge 3$) the alternating group
in $n$ letters, $\mu_I = \langle \text{\rm exp}(2 \pi \sqrt{-1})/I \rangle$
the multiplicative group of order $I$ in ${\bold C}^*$
and $C_n$ an abstract cyclic group of order $n$.

\par \noindent
{\bf 2.} For a group $G$, we write $G = A . B$ if $A$ is normal in $G$
so that $G/A = B$. We write $G = A : B$ if assume further that
$A$ is normal in $G$ and $B$ is a subgroup of $G$ such that
the composition $B \rightarrow G \rightarrow G/A = B$ is the identity
(we say then that $G$ is a {\bf semi-direct product} of $A$ and $B$).

\par \noindent
{\bf 3.} For groups $H \le G$ and $x \in G$ we denote by
$c_x : H \rightarrow G$ ($h \mapsto c_x(h) = x^{-1} h x$)
the {\bf conjugation} map.

\par \noindent
{\bf 4.} For a $K3$ surface $X$, we let $S_X$ and $T_X$ be
the Neron-Severi and transcendental lattices. 
So the $K3$ lattice $H^2(X, {\bold Z})$
contains $S_X \oplus T_X$ as a sublattice of finite index.

\par \vskip 1pc \noindent
{\bf Acknowledgement.} This work was done during
the author's visit to Hokkaido University, University of Tokyo
and Korea Institute for Advanced Study in the summer of 2004. The author would like to
thank the institutes and Professors F. Catanese, Alfred Chen, 
I. Dolgachev, J. Keum, S. Kondo,
M. L. Lang, K. Oguiso and I. Shimada for the hospitality and valuable suggestions.

\par \vskip 2pc \noindent
{\bf \S 1. Preliminary Results}

\par \vskip 1pc \noindent
{\bf (1.0).} In this section, we will prepare some basic results to be used
late. Let $X$ be a $K3$ surface with a non-zero 2-form $\omega_X$
and let $G \subseteq \text{Aut}(X)$ be a finite
group of automorphisms. For every $h \in G$,
we have $h^* \omega_X = \alpha(h) \omega_X$ for some scalar
$\alpha(h) \in {\bold C}^*$.
Clearly, $\alpha : G \rightarrow {\bold C}^*$
is a homomorphism. A fact in basic group theory says that $\alpha(G)$ is
a finite cyclic group, so 
$\alpha(G) = \mu_I = \langle \text{\rm exp}(2 \pi \sqrt{-1}/I) \rangle$ for some $I \ge 1$.
This natural number $I = I(G)$ is called the {\bf transcendental}
value of $G$. It is known that $I = I(G)$ for some $G$
if and only if that the Euler function $\varphi(I) \le 21$
and $I \ne 60$ [MO].

\par \noindent
Set $G_N = \text{\rm Ker}(\alpha)$. Then we have the {\bf basic exact sequence} below:
$$1 \longrightarrow G_N \longrightarrow G 
\overset{\alpha}\to{\longrightarrow} \mu_I \longrightarrow 1.$$

\par \noindent
For the $G$ in the basic exact sequence, we write $G = G_N . \mu_I$, though
there is no guarantee that $G = G_N : \mu_I$
(a semi-direct product).

\par \vskip 1pc \noindent
{\bf 1.0A.} If $G$ is a finite perfect group, i.e., the commutator
group $[G, G] = G$ (especially, if $G$ is a non-abelian simple group
like $A_5$), then $G = G_N$. 

\par \noindent
{\bf 1.0B.} $G_N$ acts trivially on the transcendental lattice $T_X$
(Lefschetz theorem on $(1, 1)$-classes).

\par \noindent
{\bf 1.0C.} If a subgroup $H \le G_N$ fixes a point $P$, then
$H < SL(T_{X, P}) \cong SL_2({\bold C})$ [Mu1, {\bf (1.5)}]. 
The finite subgroups of
$SL_2({\bold C})$ are listed up in [Mu1, {\bf (1.6)}]. These are
cyclic $C_n$, binary dihedral (or quaternion) $Q_{4n}$ ($n \ge 2$), binary tetrahedral $T_{24}$, 
binary octahedral $O_{48}$ and binary icosahedral $I_{120}$.

\par \vskip 1pc \noindent
{\bf Lemma 1.1.} Suppose that $G := A_5 . \mu_I$ acts faithfully
on a $K3$ surface $X$.

\par \noindent
(1) The Picard number $\rho(X) \ge 19$, and
$I = 1, 2, 3, 4, 6$. Moreover, $\rho(X) = 20$ if $I \ge 3$.

\par \noindent
(2) We have $G = A_5 : \mu_I$, i.e.,
a semi-prodcut of a normal subgroup $A_5$ and a subgroup $\mu_I$ of $G$.
Moreover, $G = A_5 \times \mu_I$ if $I = 3$.

\par \vskip 1pc \noindent
{\it Proof.} (1) In notation of [Xi, the list], 
$\rho(X) =  \text{\rm rank} \, S_X \ge c+1 = 19$.
Also the Euler function $\varphi(I)$ divides $\text{\rm rank} \, T_X = 22 - \rho(X)$ 
by [Ni1, Theorem 0.1]. So (1) follows.

\par \noindent
(2) Let $g \in G$ such that $\alpha(g)$ is a generator of $\mu_I$.
Since $\text{Aut}(A_5) = S_5 > A_5$ and the conjugation homomorphism 
$A_5 \rightarrow \text{Aut}(A_5)$
($x \mapsto c_x$) is an isomorphism onto $A_5$, 
the conjugation map $c_g$ equals $c_{(12)a}$ or $c_{a}$ on $A_5$
for some $a \in A$. Replacing $g$ by $g a^{-1}$, we may assume that
$c_g = c_{(12)}$ or $c_{\text{\rm id}}$. Thus $g^2$ commutes with every element
in $A_5$. If $2 | I$, then $g^I \in \text{\rm Ker}(\alpha) = A_5$
is in the centre of $A_5$ (which is trivial) and hence $\text{\rm ord} (g) = I$;
thus $G = A_5 : \mu_I$. If $I = 3$, then $\text{gcd}(3, \text{\rm ord} (g)/3) = 1$
as proved in [IOZ] or [Og, Proposition 5.1]; so replacing $g$
by $g^{\ell}$ with $\ell = \text{\rm ord} (g)/3$ (or $2 \text{\rm ord} (g)/3$),
we have $G = A_5 \times \langle g \rangle = A_5 \times \mu_3$.

\par \vskip 1pc \noindent
The second result below [Ni1, \S 5] 
is crucial in classifying symplectic groups in [Mu1].
For the first, see [Ni2], [Z1] or [Z2, Lemma 1.2],
where the Hodge index theorem is also used here.

\par \vskip 1pc \noindent
{\bf Lemma 1.2.} (1) Let $h$ be a non-symplectic involution on a $K3$ surface $X$.
Then $X^h$ is a disjoint union of $s$ smooth curves $C_i$
with $0 \le s \le 10$.
To be precise, $X^h$ (if not empty) is either a disjoint union of a genus $\ge 2$
curve $C$ and a few ${\bold P}^1$'s, or
a disjoint union of a few elliptic curves and ${\bold P}^1$'s,
or a disjoint union of a few ${\bold P}^1$'s.

\par \noindent
(2) If $\delta$ is a non-trivial symplectic automorphism of finite 
order on a $K3$ surface $X$, then $\text{\rm ord} (\delta) \le 8$ and
$X^{\delta}$ is a finite set. 
To be precise, if $\text{\rm ord} (\delta) = 2, 3, 4, 5, 6, 7, 8$,
then $|X^{\delta}| = 8, 6, 4, 4, 2, 3, 2$, respectively.
In particular, if $A_5 \subseteq \text{Aut}(X)$ then
$\sum_{\delta \in A_5} \chi_{\text{\rm top}}(X^{\delta}) = 360$ (see {\bf (1.0A)}).

\par \vskip 1pc \noindent
For an automorphism $h$ on a smooth algebraic surface $Y$, we split the
pointwise fixed locus as the disjoint union
of 1-dimensional part and the isolated part:
$Y^h = Y^h_{1-\text{\rm dim}} \coprod Y^h_{\text{\rm isol}}$.
The proof of (1) below is similar to that for (1) in {\bf (1.2)}.

\par \vskip 1pc \noindent
{\bf 1.3.} (1) $Y^h_{1-\text{\rm dim}}$ (if not empty) is a disjoint union of smooth curves.

\par \noindent
(2) The Euler number $\chi_{\text{\rm top}}(Y^h_{1-\text{\rm dim}}) = \sum_{C}
(2 - 2 g(C)) = 2n_h$ for some integer $n_h$, where
$C$ runs in the set $Y^h_{1-\text{\rm dim}}$ of curves.

\par \noindent
(3) The Euler number $\chi_{\text{\rm top}}(Y^h) = m_h + 2 n_h$, where
$m_h = |Y^h_{\text{\rm isol}}|$.

\par \vskip 1pc \noindent
The results of [IOZ] below follow from the application of
Lefschetz fixed point 
formula to the trivial vector bundle in Atiyah-Segal-Singer
[AS2, AS3, pages 542 and 567]. 
For a proof, see [OZ1, Lemma 2.3] and [Z2, Proposition 1.4].

\par \vskip 1pc \noindent
{\bf Lemma 1.4.} Let $X$ be a $K3$ surface and let
$h \in \text{Aut}(X)$ be of order $I$ such that $h^* \omega_X = \eta_I \omega_X$
for some primitive $I$-th root $\eta_I$ of 1.

\par \noindent
(1) Suppose that $I = 3$. Then $m_h = 3 + n_h$ and hence 
$\chi_{\text{\rm top}}(X^h) = 3(1 + n_h)$. Moreover, $-3 \le n_h \le 6$.

\par \noindent
(2) Suppose that $I = 3$. 
If $\delta \in \text{Aut}(X)$ is symplectic of order 5 and
commutes with $h$. Then $|X^{h \delta}| = 4$.

\par \vskip 1pc \noindent
The following result can be found in [Ni1, Theorem 0.1], [MO, Lemma (1.1)],
or [OZ3, Lemma (2.8)].

\par \vskip 1pc \noindent
{\bf Lemma 1.5.} Suppose that $X$ is a $K3$ surface of Picard
number $\rho(X) = 20$ and $g$ an order-3 automorphism
such that $g^* \omega_X = \eta_3 \omega_X$ with a primitive
3rd root $\eta_3$ of 1. Then we can express the transcendental lattice
$T_X$ as $T_X = {\bold Z}[t_1, t_2]$ so that
$t_2 = g^*(t_1)$, $g^*(t_2) = -(t_1+t_2)$. In particular,
for some $m \ge 1$, 
the intersection form $(t_i . t_j) = \pmatrix 2m & -m \\ -m & 2m \endpmatrix$.

\par \vskip 1pc \noindent
Now we assume that $G = G_N . \mu_I$ (with $I = I(G)$)
acts on a $K3$ surface $X$.
When $G_N = A_5$,
we will determine the action of $G_N$
on the Neron Severi lattice $S_X$ of $X$:

\par \vskip 1pc \noindent
{\bf Lemma 1.6.} (1) Suppose that $A_5$ acts on a $K3$ surface $X$, 
and $\text{\rm rank} \, S_X = 20$ 
(this is true if $I \ge 3$
by {\bf (1.1)}). Then
we have the irreducible decomposition below (in the
notation of Atlas for the characters of $A_5$),
where $\chi_1$ (the trivial character), 
$\chi_4$ and $\chi_5$ have dimensions 1, 4 and 5, respectively,
where $\chi_i'$ is a copy of $\chi_i$:
$$S_X \otimes {\bold C}
= \chi_1 \oplus \chi_1' \oplus \chi_4 \oplus \chi_4' \oplus \chi_5 \oplus \chi_5'.$$

\par \noindent
(2) For conjugacy class $nA$ (and $nB$) of order $n$ in $A_5$ and the characters 
$\chi_i$ of $A_5$, we have the following {\bf Table 1} from [Atlas], where $Z$ is respectively 
$1A$, $2A$, $3A$, $5A$ or $5B$:
$$\gather
[\chi_1, \chi_2, \chi_3, \chi_4, \chi_5](Z) =
[1, 3, 3, 4, 5], \hskip 1pc [1, -1, -1, 0, 1], \hskip 1pc [1, 0, 0, 1, -1], \\
[1, \, (1- \sqrt{5})/2, \,\,(1 + \sqrt{5})/2, \, -1, 0], \hskip 1pc
[1, \, (1+ \sqrt{5})/2, \,\,(1 - \sqrt{5})/2, \, -1, 0].
\endgather $$

\par \vskip 1pc \noindent
{\it Proof.} The assertion(1) appeared in [Z2].
For the readers' convenience, we reprove it here.
Applying the Lefschetz fixed point formula to the action
of $A_5$ on 
$H^*(X, {\bold Z}) = \oplus_{i=0}^4 H^i(X, {\bold Z})$
and noting that $H^2(X, {\bold Z})$
contains $S_X \oplus T_X$ as a finite index sublattice, we obtain
(see also {\bf (1.0A-B)} and {\bf (1.2)}):
$$2 + \text{\rm rank} \, T_X + \text{\rm rank} (S_X)^{A_5} = 
\text{\rm rank} \,H^*(X, {\bold Z})^{A_5} = \frac{1}{|A_5|} \sum_{a \in A_5} \chi_{\text{\rm top}}(X^a) = 360/60 = 6.$$

\par \noindent
Thus $\text{\rm rank} \, S_X^{A_5} = 2$.
So the irreducible decomposition is of the
following form, where $a_{i}$ are non-negative integers:
$$S(X) \otimes {\bold C} = 2 \chi_{1} \oplus a_{2}\chi_{2} 
\oplus a_{3}\chi_{3} \oplus a_{4}\chi_{4} \oplus a_{5}\chi_{5}.$$
Using the topological Lefschetz
fixed point formula, the fact that $\text{\rm rank}\, T(X)$
$= 2$ and {\bf (1.0B)}, we have, for $a \in A_{5}$, that:
$$\chi_{\text{\rm top}}(X^{a}) = 2 + \text{\rm rank} \, T_X + 
\text{\rm Tr}(a^{*} \vert S(X))$$
Running $a$ through the five conjugacy classes and
calculating both sides, using {\bf (1.2)} and the character Table 1 in (2),
we obtain the following system of equations:
$$\gather 20 = 2 + 3(a_{2} + a_{3}) + 4a_{4} + 5a_{5}, \\
4 = 2 - (a_{2} + a_{3}) + a_{5}, \\
2 = 2 +  a_4 - a_{5}, \\
0 = 2 + \frac{1 - \sqrt{5}}{2}a_{2} + \frac{1 + \sqrt{5}}{2} a_{3} 
- a_{4}, \\
0 = 2 +  \frac{1 + \sqrt{5}}{2}a_{2} + \frac{1 - \sqrt{5}}{2} a_{3} 
- a_{4}.
\endgather$$

\par \noindent
Now, we get the result by solving this system of Diophantine equations.

\par \vskip 1pc \noindent
The two results below are either easy or well known and will be
frequently used in the arguments of the subsequent sections.

\par \vskip 1pc \noindent
{\bf Lemma 1.7.} Let $f : A_5 \rightarrow S_r$ ($r \ge 2$) be a
homomorphism.

\par \noindent
(1) If $r = 2$, $3$, or $4$, then $f$ is trivial.

\par \noindent
(2) If $\text{\rm Im}(f)$ is a transitive subgroup of the full symmetry group $S_r$
in $r$ letters $\{1, 2, \dots, r\}$ (whence $r \ge 5$ by (1)), then
$r | |A_5|$ with $|A_5|/r$ equal to the order of the subgroup of 
$A_5$ stabilizing the letter 1, so $r \in \{5, 6, 10, 12, 15, 20, 30\}$. 

\par \vskip 1pc \noindent
{\bf Lemma 1.8.}  
(1) $\text{Aut}({\bold P}^1)$ includes $A_5$ but not $S_5$ [Su, Theorem 6.17].
The action of $A_5$ on ${\bold P}^1$ is unique up to isomorphisms.

\par \noindent
(2) Every $A_5$ in $PGL_3({\bold C})$ is the image of
an $A_5$ in $SL_3({\bold C})$.

\par \noindent
(3) The conjugation by $(12) \in S_5$ switches the two 3-dimensional characters 
$\chi_2$ and $\chi_3$ of $A_5$ [Atlas].

\par \noindent
(4) If $\text{\rm id} \ne f \in \text{Aut}({\bold P}^1)$ is an automorphism of finite order,
then $f$ fixes exactly two point of ${\bold P}^1$
(by the diagonalization of a lifting of $f$ to $GL_2({\bold C})$).

\par \noindent
(5) If $f_r$ ($r = 2$ or $3$) is an order$-r$ automorphism of an 
elliptic curve $E$, then either $f_r$ acts freely on $E$, or
the fix locus satisfies
$|X^{f_r}| = 4$ (resp. $= 3$) if $r = 2$ (resp. $r = 3$) (by the Hurwitz formula).

\par \vskip 1pc \noindent
{\it Proof.} (1) For the uniqueness of the action of $A_5$ on ${\bold P}^1$,
one may assume the representation of 
$D_{10} = \langle \gamma = (12345), \sigma = (14)(23) \rangle$
is given by $\gamma : z \rightarrow \eta z$ with $\eta$ a primitive 5-th root of 1
and $\sigma : z \rightarrow \alpha/z$. Note that $A_5 = \langle \gamma, \varepsilon \rangle$
with $\varepsilon = (12)(34)$. If one lets $\varepsilon : z \rightarrow (az + b)/(c z + d)$
be in $\text{Aut}({\bold P}^1)$, then one can check that 
$d = -a$ because $\text{ord}(\varepsilon) = 2$, and also
$b = - c \alpha$ because $\varepsilon$ commutes with $\sigma$.
So $\varepsilon : z \rightarrow (z - \alpha e)/(e z - 1)$ with $e = c/a$.
The commutativity of $\varepsilon \sigma \gamma^2 \varepsilon = (12)(45)$
with $\sigma \gamma^{-1} = (15)(24)$ implies that $e^2 \alpha = \eta + \eta^{-1} - 1$.
Now let $\rho : z \rightarrow e \alpha/z$ be in $\text{Aut}({\bold P}^1)$.
Then $\rho^{-1} \gamma \rho : z \rightarrow \eta^{-1} z$,
$\rho^{-1} \sigma \rho : z \rightarrow e^2 \alpha/z$
and $\rho^{-1} \varepsilon \rho : z \rightarrow (z - e^2 \alpha)/(z-1)$.
Hence the action of $A_5$ on ${\bold P}^1$ is unique modulo isomorphisms.

\par \noindent
(2) For an $A_5$ in $SL_3({\bold C})$, see [Bu, \S 232].
The inverse $\widetilde{A}_5 \subset SL_3({\bold C})$
of an $A_5 \subset PGL_3({\bold C})$ is of the form
$\widetilde{A}_5 = A_5 : \mu_3$ (indeed, a direct product) because
the Schur multiplier $M(A_5) = 2$, coprime to $3$ [Atlas]. So (2) follows.

\par \vskip 2pc \noindent
{\bf \S 2. Alternating groups actions on the Niemeier lattices}

\par \vskip 1pc \noindent
For a $K3$ surface $X$, denote 
by $L = H^2(X, {\bold Z})$
the {\bf K3 lattice},
$S_X = \text{Pic} \, X$ (now) the {\bf Neron-Severi lattice} and $T_X$ the 
{\bf transcendental lattice}.
So $T_X = S_X^{\perp}$ in $L$ and $L$ contains a finite-index sublattice
$S_X \oplus T_X$.

\par \vskip 1pc \noindent
{\bf (2.0).} Suppose that $G_N = A_5$ acts faithfully on $X$.
In this section we shall prove Theorem C which is part of {\bf (2.1)} below.
Indeed, by the proof of {\bf (1.6)}, we have
$\text{\rm rank} \, L^{G_N} =
\text{\rm rank} \, T_X + \text{\rm rank} \, S_X^{G_N} = 4$,
so $(\text{\rm rank} \, T_X, \, \text{\rm rank} \, S_X, \, \text{\rm rank} \, S_X^{G_N})
= (2, 20, 2)$ or $(3, 19, 1)$.

\par \vskip 1pc \noindent
Denote by $L^{G_N} := \{x \in L \,|\, g^*x = x$
for all $g \in G_N\}$ and its orthogonal $L_{G_N} := (L^{G_N})^{\perp} =
\{x \in L \,|\, (x, y) = 0$ for all $y \in L^{G_N} \}$. 
Then $L^{G_N}$ contains $S_X^{G_N} \oplus T_X$ as a sublattice of finite index
by {\bf (1.0A-B)}. 

\par \noindent 
By [Ko1, Lemmas 5 and 6], there are a (non-Leech) Niemeier lattice $N(Rt)$,
a primitive embedding $A_1 \oplus L_{G_N} \subset N(Rt)$ 
and a faithful action of $G_N$ on $N(Rt)$
such that $L_{G_N} = N(Rt)_{G_N}$, and
the action of $G_N$ on the summand $A_1$ is trivial
and stabilizes a Weyl chamber (one of whose codimension one faces
corresponds to this $A_1$). Moreover,
$G_N \le S(N(Rt)): = O(N(Rt))/W(N(Rt)) \le O(Rt)/W(N(Rt)) (=: \text{\rm Sym}(Rt))$,
where $\text{\rm Sym}(Rt)$ is the full symmetry group of the Coxeter-Dynkin diagram $Rt$.
Note that $\text{\rm rank} \, N(Rt)^{G_N} =
2 + \text{\rm rank} \, L^{G_N} = 6$
and the discriminant groups satisfy:
$$(*) \hskip 1pc
A_{L^{G_N}} \cong A_{L_{G_N}} (-1) = A_{N(Rt)_{G_N}} (-1) \cong A_{N(Rt)^{G_N}}.$$

\par \vskip 1pc \noindent
Now
$N(Rt)^{G_N}$ is a \text{\rm rank} $6$ lattice generated by $e_1, \dots, e_6$ say.
Denote by $M = (e_i . e_j)$ the intersection matrix and $M^{-1} =
(f_1, \dots, f_6)$ with $f_j$ column vectors and set 
$e_i^* = (e_1, \dots, e_6) f_i$. Then 
$(N(Rt)^{G_N})^{\vee} = \text{\rm Hom}(N(Rt)^{G_N}, {\bold Z})$ has the dual basis
$\{e_1^*, \dots, e_6^*\}$ with the intersection matrix
$(e_i^* . e_j^*)_{1 \le i, j \le 6} = M^{-1}$. 
The discriminant groups 
satisfy $A_{L^{G_N}} \cong
A_{N(Rt)^{G_N}} = {\bold Z}[e_1^*, \dots, e_6^*]/{\bold Z}[e_1, \dots, e_6]$.

\par \vskip 1pc \noindent
In this section, we shall prove the following result
(much inspired by [Ko1]), which (and the proof of which)
should be useful in studying $\text{Aut}(X)$ from the $K3$ lattice point of view.
This result is used in [Z2, Lemma 3.5].

\par \vskip 1pc \noindent
{\bf Theorem 2.1.} Suppose that $G_N = A_5$ acts faithfully on a $K3$ surface $X$. 

\par \noindent
(1) We have $Rt = 24A_1$ or $Rt = 12 A_2$.
The lattice $N(Rt)^{G_N}$ is of rank $6$ and generated by $e_1, \dots, e_6$ say.
Denote by $M = (e_i . e_j)$ the intersection matrix and write $M^{-1} =
(f_1, \dots, f_6)$ with $f_j$ column vectors and set $e_i^* = (e_1, \dots, e_6) f_i$.

\par \noindent
(2) If $Rt = 24A_1$, then the orbit decomposition of the $G_N$-action on the
$24$ simple roots is either one of 
$$\text{\rm (i)} \,\, [1, 1, 5, 5, 6, 6], \,\,\,\, \text{\rm (ii)} \,\, [1, 1, 1, 5, 6, 10], \,\,\,\,
\text{\rm (iii)} \,\, [1, 1, 1, 1, 5, 15], \,\,\,\, \text{\rm (iv)} \,\, [1, 1, 1, 1, 10, 10].$$

\par \noindent
If $Rt = 12A_2$, then the orbit decomposition of the $G_N$-action on
the $24$ simple roots is either one of
$$(v) \,\, [1, 1, 1, 1, 10, 10], \hskip 1pc (vi) \,\, [1,1,5,5,6,6],$$
where in (v) (resp. (vi)) $10A_2$ (resp. $5A_2$, or $6A_2$)
is split into two orbits with 10 (resp. 5, or 6) disjoint roots each.

\par \noindent
(3) For Case(2i), the intersection matrix $M_1 = (e_i . e_j)$ and its inverse 
$M_1^{-1}$ are respectively:
$$\pmatrix
-2 & 0 & 0 & -1 & -1 & -1 \\
0 & -2 & 0 & -1 & -1 & -1 \\
0 & 0 & -10 & 0 & 0 & -5 \\
-1 & -1 & 0 & -4 & -1 & -1 \\
-1 & -1 & 0 & -1 & -4 & -1 \\
-1 & -1 & -5 & -1 & -1 & -6
\endpmatrix,
\pmatrix
-23/30 & -4/15 & -1/10 & 1/6 & 1/6 & 1/5 \\
-4/15 & -23/30 & -1/10 & 1/6 & 1/6 & 1/5 \\
-1/10 & -1/10 & -1/5 & 0 & 0 & 1/5 \\ 
1/6 & 1/6 & 0 & -1/3 & 0 & 0 \\
1/6 & 1/6 & 0 & 0 & -1/3 & 0 \\
1/5 & 1/5 & 1/5 & 0 & 0 & -2/5 \endpmatrix.
$$
The discriminant group (cf. {\bf (2.0)}, $A_{L^{G_N}} \cong$) 
$A_{N(Rt)^{G_N}} = 
\text{\rm Hom}(N(Rt)^{G_N}, \, {\bold Z})/N(Rt)^{G_N}$ 
$\cong {\bold Z}/(30) \oplus {\bold Z}/(30)$
and is generated by cosets $\overline{e}_1^*$ and 
$\overline{e}_2^* + \overline{e}_3^* + \overline{e}_4^*$ with intersection form:
$$\pmatrix (\overline{e}_1^*)^2 & 
\overline{e}_1^* . (\overline{e}_2^* + \overline{e}_3^* + \overline{e}_4^*) \\
\overline{e}_1^* . (\overline{e}_2^* + \overline{e}_3^* + \overline{e}_4^*) &
(\overline{e}_2^* + \overline{e}_3^* + \overline{e}_4^*)^2 \endpmatrix =
\pmatrix
-23/30 & -1/5 \\ -1/5 & -35/30
\endpmatrix.$$

\par \noindent
(4) For Case(2ii), the intersection matrix $M_2 = (e_i . e_j)$ and its inverse
$M_2^{-1}$ are respectively:
$$\pmatrix
-2 & 0 & 0 & -1 & -1 & -1 \\
0 & -2 & 0 & -1 & -1 & -1 \\
0 & 0 & -2 & -1 & 0 & 0 \\
-1 & -1 & -1 & -4 & -1 & -1 \\
-1 & -1 & 0 & -1 & -4 & -1 \\
-1 & -1 & 0 & -1 & -1 & -6
\endpmatrix,
\pmatrix
-11/15 & -7/30 & -1/10 & 1/5 & 1/6 & 1/10 \\
-7/30 & -11/15 & -1/10 & 1/5 & 1/6 & 1/10 \\
-1/10 & -1/10 & -3/5 & 1/5 & 0 & 0 \\
1/5 & 1/5 & 1/5 & -2/5 & 0 & 0 \\
1/6 & 1/6 & 0 & 0 & -1/3 & 0 \\
1/10 & 1/10 & 0 & 0 & 0 & -1/5
\endpmatrix.
$$
The discriminant group $A_{N(Rt)^{G_N}}$ is isomorphic to 
${\bold Z}/(30) \oplus {\bold Z}/(10)$
and generated by the cosets 
$\overline{e}_1^*$ and $\overline{e}_3^*$
with intersection form:
$$\pmatrix (\overline{e}_1^*)^2 & \overline{e}_1^* . \overline{e}_3^* \\ 
\overline{e}_1^* . \overline{e}_3^* & (\overline{e}_3^*)^2 \endpmatrix =
\pmatrix -11/15 & -1/10 \\ -1/10 & -3/5 \endpmatrix.$$

\par \noindent
(5) For Case(2iii), the intersection matrix $M_3 = (e_i . e_j)$ and
$M_3^{-1}$ are respectively:
$$\pmatrix
-2 & 0 & 0 & 0 & -1 & 0 \\
0 & -2 & 0 & 0 & -1 & 0 \\
0 & 0 & -2 & 0 & -1 & 0 \\ 
0 & 0 & 0 & -2 & 0 & -1 \\
-1 & -1 & -1 & 0 & -4 & 0 \\
0 & 0 & 0 & -1 & 0 & -8
\endpmatrix,
\pmatrix
-3/5 & -1/10 & -1/10 & 0 & 1/5 & 0 \\ 
-1/10 & -3/5 & -1/10 & 0 & 1/5 & 0 \\ 
-1/10 & -1/10 & -3/5 & 0 & 1/5 & 0 \\
0 & 0 & 0 & -8/15 & 0 & 1/15 \\
1/5 & 1/5 & 1/5 & 0 & -2/5 & 0 \\
0 & 0 & 0 & 1/15 & 0 & -2/15
\endpmatrix. $$
The discriminant group $A_{N(Rt)^{G_N}}$ is isomorphic to 
${\bold Z}/(30) \oplus {\bold Z}/(10)$
and generated by the cosets
$\overline{e}_2^*$ and $\overline{e}_1^* + \overline{e}_4^*$
with intersection form:
$$\pmatrix (\overline{e}_2^*)^2 & 
\overline{e}_2^* . (\overline{e}_1^* + \overline{e}_4^*)\\ 
\overline{e}_2^* . (\overline{e}_1^* + \overline{e}_4^*) & 
(\overline{e}_1^* + \overline{e}_4^*)^2 \endpmatrix =
\pmatrix -3/5 & -1/10 \\ -1/10 & 13/15 \endpmatrix.$$

\par \noindent
(6) For Case(2iv), the intersection matrix $M_4 = (e_i . e_j)$ and
$M_4^{-1}$ are respectively:
$$\pmatrix
-2 & 0 & 0 & 0 & -1 & 0 \\
0 & -2 & 0 & 0 & -1 & 0 \\
0 & 0 & -2 & 0 & 0 & -1 \\
0 & 0 & 0 & -2 & 0 & -1 \\
-1 & -1 & 0 & 0 & -6 & 0 \\
0 & 0 & -1 & -1 & 0 & -6
\endpmatrix,
\pmatrix
-11/20 & -1/20 & 0 & 0 & 1/10 & 0 \\
-1/20 & -11/20 & 0 & 0 & 1/10 & 0] \\
0 & 0 & -11/20 & -1/20 & 0 & 1/10 \\
0 & 0 & -1/20 & -11/20 & 0 & 1/10 \\
1/10 & 1/10 & 0 & 0 & -1/5 & 0 \\
0 & 0 & 1/10 & 1/10 & 0 & -1/5
\endpmatrix.$$
The discriminant group $A_{N(Rt)^{G_N}}$ is isomorphic to 
${\bold Z}/(20) \oplus {\bold Z}/(20)$
and generated by the cosets 
$\overline{e}_1^*$ and $\overline{e}_3^*$
with intersection form:
$$\pmatrix (\overline{e}_1^*)^2 & 
\overline{e}_1^* . \overline{e}_3^* \\ 
\overline{e}_1^* . \overline{e}_3^* & 
(\overline{e}_3^*)^2 \endpmatrix =
\pmatrix -11/20 & 0 \\ 0 & -11/20 \endpmatrix.$$

\par \noindent
(7) For Case(2v), the intersection matrix $M_5 = (e_i . e_j)$ and its inverse $M_5^{-1}$
are respectively:
$$\pmatrix
-2 &1 &0 &0 &0 &0 \\ 
1 &-2 &0 &0 &0 &-1 \\ 
0 &0 &-2 &1 &0 &0 \\ 
0 &0 &1 &-2 &0 &-1 \\
0 &0 &0 &0 &-20 &0 \\ 
0 &-1 &0 &-1 &0 &-8
\endpmatrix,
\pmatrix
-41/60 &  -11/30 &  -1/60 &  -1/30 &  0 &  1/20 \\
-11/30 &  -11/15 &  -1/30 &  -1/15 &  0 &  1/10 \\ 
-1/60 &  -1/30 &  -41/60 &  -11/30 &  0 &  1/20 \\ 
-1/30 &  -1/15 &  -11/30 &  -11/15 &  0 &  1/10 \\
0 &  0 &  0 &  0 &  -1/20 &  0 \\
1/20 &  1/10 &  1/20 &  1/10 &  0 &  -3/20
\endpmatrix.
$$
The discriminant group $A_{N(Rt)^{G_N}}$ is isomorphic to 
${\bold Z}/(60) \oplus {\bold Z}/(20)$
and generated by the cosets 
$\overline{e}_1^*$ and $\overline{e}_5^*$
with intersection form:
$$\pmatrix (\overline{e}_1^*)^2 & \overline{e}_1^* . \overline{e}_5^* \\ 
\overline{e}_1^* . \overline{e}_5^* & (\overline{e}_5^*)^2 \endpmatrix =
\pmatrix -41/60 & 0 \\ 0 & -1/20 \endpmatrix.$$

\par \noindent
(8) For Case(2vi), The intersection matrix $M_6 = (e_i . e_j)$ and $M_6^{-1}$
are respectively:
$$\pmatrix
-2 & 1 & 0 & 0 & 0 & 0 \\
1 & -2 & 0 & 0 & -1 & 0 \\ 
0 & 0 & -10 & 0 & 0 & 0 \\ 
0 & 0 & 0 & -12 & 0 & 0 \\ 
0 & -1 & 0 & 0 & -4 & 0 \\ 
0 & 0 & 0 & 0 & 0 & -4
\endpmatrix,
\pmatrix
-7/10 & -2/5 & 0 & 0 & 1/10 & 0 \\ 
-2/5 & -4/5 & 0 & 0 & 1/5 & 0 \\ 
0 & 0 & -1/10 & 0 & 0 & 0 \\
0 & 0 & 0 & -1/12 & 0 & 0 \\
1/10 & 1/5 & 0 & 0 & -3/10 & 0 \\
0 & 0 & 0 & 0 & 0 & -1/4
\endpmatrix.
$$
The discriminant group $A_{N(Rt)^{G_N}} = {\bold Z}/(60) \oplus {\bold Z}/(20)
\oplus {\bold Z}/(2) \oplus {\bold Z}/(2) = {\bold Z}/(10) \oplus {\bold Z}/(10)
\oplus {\bold Z}/(12) \oplus {\bold Z}/(4)$
and the latter is generated by the cosets $\overline{e}_j^*$ ($j = 1, 3, 4, 6$).

\par \noindent
(9) In both of the cases of $M_2$ and $M_3$, the discriminant
group $A_{N(Rt)^{G_N}}$ is isomorphic to the group
$\langle \overline{t}_1^*, \overline{t}_2^* \rangle \cong 
{\bold Z}/(30) \oplus {\bold Z}/(10)$ with the intersection matrix
$(\overline{t}_i^* . \overline{t}_j^*) = \pmatrix
1/15 & 1/30 \\ 1/30 & 1/15 \endpmatrix$.

\par \vskip 1pc \noindent
We now prove {\bf (2.1)}.
Since $\text{\rm rank} \, N(Rt)^{G_N} = 6$, the $G_N$-action
on the $24$ simple roots of $Rt$ has exactly 6 orbits.

\par \noindent 
We argue as in the proof of
[Ko1, Theorem 4]. The fact that $G_N = A_5 < S(N(Rt))$
implies that $Rt$ is one of the following:
$24A_1, 12A_2, 6A_4, 6D_4$.

\par \noindent
If $Rt = 6A_4$, then $S(N(Rt)) = 2. PGL_2(5)$ ($< 2 . S_6$) [CS, Ch 16, \S 1],
where the order 2 element acts as a symmetry of order 2 on each connected component of
Dynkin type $A_4$, and $PGL_2(5)$ acts on the set (identified with
$\{0, 1, 2, 3, 4, \infty\}$) of 6 components of $Rt$
as permutations in a natural way. Since
$A_5$ is simple, the composition of homomorphisms below is an injection:
$A_5 \subset S(N(Rt)) \rightarrow PGL_2(5)$,
so we may assume that $A_5 < PGL_2(5)$. Since $G_N = A_5$ fixes one simple
root of $Rt$ by the construction, our $A_5$ is a subgroup of the stabilizer
subgroup of $PGL_2(5)$ and this stablizer is of order $|PGL_2(5)|/6 = 20$.
This is impossible because $|A_5| = 60 > 20$.

\par \noindent
If $Rt = 6D_4$, then $S(N(Rt)) = 3 . S_6$, where the 
order 3 element acts as a symmetry of order 3 on each connected component of
Dynkin type $D_4$, and $S_6$ acts on the set of 6 connected components of $Rt$
as permutations. As above, the simplicity of $G_N$ implies that the subgroup
$G_N$ of $S(N(Rt))$ is indeed a subgroup of $S_6$.
Since $G_N = A_5$ fixes one simple root of $Rt$, our group $A_5$
is a subgroup ($=[S_5, S_5]$) of the stabilizer subgroup $S_5$ of $S_6$.
So this $A_5$ acts transitively on
the remaining 5 connected components of $Rt$ and hence the $G_N$-action on the 24
simple roots has exactly 8 orbits, noting that
one connected component of $Rt$ is component wise fixed by $G_N$, a contradiction.

\par \noindent
Suppose that $Rt = 12A_2$. Then $S(N(Rt)) = 2. M_{12}$, where the 
order 2 element acts as a symmetry of order 2 on each connected component of
Dynkin type $A_2$, and the Mathieu group $M_{12}$ acts on the set of 12 
connected components of $Rt$
as permutations. Let $r_{2k-1} + r_{2k}$ ($1 \le k \le 12$) be
the 12 connected components of $Rt$ with $r_j$ the 24 simple roots.
Every non-trivial element of $N(Rt)/Rt$ is of the 
form $\sum_{i \in H} \pm (r_{2i-1} + 2r_{2i})/3$ where $H$ is an element
of the ternary Golay code and $|H| = 6, 9, 12$ [CS, Ch 3, \S 2.8.5].
Since the group $G_N = A_5$ is simple and fixes one simple root of $Rt$, this $G_N$
is a subgroup of $M_{12}$ and indeed, a subgroup of the stabilizer subgroup $M_{11}$
of $M_{12}$. Suppose the $G_N$-orbit decomposition on the $12$ connected components
is $1 + a + b$. Then the $G_N$-orbit decomposition of the 24 simple roots is
$[1, 1, a, a, b, b]$ (so $a + b = 11$), where $aA_2$ (resp. $bA_2$) is split into 
two $G_N$-orbits with $a$ (resp. $b$) disjoint simple roots each.
Thus Case (2v) or (2vi) occurs by {\bf (1.7)}.

\par \noindent
Suppose that $Rt = 24A_1$. Then $S(N(Rt)) = M_{24}$.
The elements of $N(Rt)/Rt$ form the binary Golay code.
Since $G_N = A_5$ fixes one simple root of $Rt$, our group $G_N$
is the stabilizer subgroup $M_{23}$ of $M_{24}$. 
Let the $G_N$-orbit decomposition of the 24 simple roots be
$[1, a, b, c, d, e]$ with $a \le b \le c \le d \le e$
(so $a + b + c + d + e = 23$).
By {\bf (1.7)}, all $a, b, c, d, e$ are in $\{1, 5, 6, 10, 12, 15, 20\}$
and hence Cases (2i) - (2vi) occur.

\par \noindent
(4) According to the ordering of $[1, 1, 1, 5, 6, 10]$, we label
the orbits as $O_1 = \{r_1\}, O_1' = \{r_2\}, O_1'' = \{r_3\}, 
O_5 = \{r_4, \dots, r_8\}, O_6 = \{r_9, \dots, r_{14}\}, O_{10} = 
\{r_{15}, \dots, r_{24}\}$, where $r_j$ are the 24 simple roots.
We claim that $O_1 + O_1' + O_1'' + O_5$ (to be precise, after divided by 2)
is an octad, and $O_1 + O_1' + O_6$ is also
an octad (after relabelling $O_1, O_1', O_1''$).
So 
$$e_i = r_i (1 \le i \le 3), \, e_4 = \frac{1}{2}(O_1 + O_1' + O_1'' + O_5), \,
e_5 = \frac{1}{2}(O_1 + O_1' + O_6), \, e_6 = \frac{1}{2}(O_1 + O_1' + O_{10})$$
form a basis of $N(Rt)^{G_N}$, noting that the last dodecad is the complement
of the symmetric sum (a dodecad) of the two octads above and that
except for the above-mentioned two octads and two dodecads,
there is no any other octad or dodecad which is a union of orbits.

\par \noindent
Indeed, let $Oct_1$ be the unique octad containing $O_5$. Note that the cycle
type in $M_{24}$ of an order-5 element $\gamma$ in $A_5$ is
$(5^4)$ , Appendix B, Table 5.I]. So $\gamma$ is of type $(5^2)$ (resp.
$(5)$) on $O_{10}$ (resp. on $O_5$ and $O_6$).
Since $\gamma(Oct_1) \cap Oct_1$ contains $O_5$, we have
$\gamma(Oct_1) = Oct_1$. If $Oct_1$ contains an element of $O_{10}$
then it contains the five images in $O_{10}$ by the action of 
$\langle \gamma \rangle$, so $|Oct_1| \ge 10$, absurd.
If $Oct_1$ contains an element $r_j$ in $O_6$ we may choose
$\gamma$ not fixing $r_j$ (note that the stabilizer subgroup of
$A_5$, regarded as a subgroup of $\text{\rm Sym}(O_6) = S_6$ and fixing
an element ($\ne r_j$) in $O_6$, has order $10$ and hence gives
rise to such $\gamma$). Then we will get a similar contradiction.
Thus $Oct_1 = O_1 + O_1' + O_1'' + O_5$ as claimed.

\par \noindent
Let $Oct_2$ be the unique octad containing the first 5 elements in $O_6$.
Let $\gamma$ be an order-5 element in $A_5$ fixing the last element in $O_6$.
Then $\gamma(Oct_2) = Oct_2$. As above, this implies that $Oct_2$ is disjoint
from $O_5$ and $O_{10}$. So either $Oct_2 = O_1 + O_1' + O_6$ after relabelling
the 1-element orbits, or $Oct_2$ is the union of the 5 elements in $O_6$
and the three 1-element orbits (this leads to that the symmetric sum
of $Oct_1$ and $Oct_2$ is a 10-word Golay code, absurd).

\par \noindent
(3) For the orbit decomposition $[1, 1, 5, 5, 6, 6]$, 
we label the orbits as $O_1 = \{r_1\}$, $O_1' = \{r_2\}$,
$O_5 = \{r_3, \dots, r_7\}$, $O_5' = \{r_8, \dots, r_{12}\}$, $O_{6} = \{r_{13}, \dots, r_{18}\}$, $O_{6}' = \{r_{19}, \dots, r_{24}\}$. As in (4), we can prove that
both $O_1 + O_1' + O_6$ and $O_1 + O_1' + O_6'$ are octads. 
Thus $N(Rt)^{G_N}$ has a basis below, noting that except for
the two octads, the symmetric sum (a dodecad) of the two octads and the
complement (another dodecad) of this dodecad, there is no other octad or
dodecad which is the union of orbits:
$$e_i = r_i (i = 1, 2), \, e_3 = O_5, \,
e_4 = \frac{1}{2}(O_1 + O_1' + O_{6}), \, e_5 = \frac{1}{2}(O_1 + O_1' + O_{6}'),
\, e_6 = \frac{1}{2}(O_1 + O_1' + O_5 + O_5')$$.

\par \noindent
(5) For the orbit decomposition $[1, 1, 1, 1, 5, 15]$, 
we label the orbits as $O_1 = \{r_1\}$, $O_1' = \{r_2\}$,
$O_1'' = \{r_3\}$, $O_1''' = \{r_4\}$, $O_{5} = \{r_5, \dots, r_{9}\}$,
$O_{15} = \{r_{10}, \dots, r_{24}\}$. As in (4), we may assume that
$O_1 + O_1' + O_1'' + O_5$ is an octad after relabelling the 1-element orbits
and that there is no any other octad or dodecad which is a union of orbits.
Thus $N(Rt)^{G_N}$ has a basis:
$$e_i = r_i (i = 1, 2, 3, 4), \,\,
e_5 = \frac{1}{2}(O_1 + O_1' + O_1'' + O_{5}), \,\, e_6 = \frac{1}{2}(O_1''' + O_{15}).$$

\par \noindent
(6) For the orbit decomposition $[1, 1, 1, 1, 10, 10]$, 
we label the orbits as $O_1 = \{r_1\}$, $O_1' = \{r_2\}$,
$O_1'' = \{r_3\}$, $O_1''' = \{r_4\}$, $O_{10} = \{r_5, \dots, r_{14}\}$,
$O_{10}' = \{r_{15}, \dots, r_{24}\}$.
Take an order-5 element $\gamma$ of $A_5$. So $O_{10}$ splits into
two 5-element subsets on each of which $\gamma$ acts transitively.
Let $Oct_j$ ($j = 1, 2$) be the unique octad containing the first (resp. second)
5-element subset. As in (4), we can show that
each $Oct_j$ is the union of the 5-element subset and three 1-element orbits.
The symmetric sum of $Oct_1$ and $Oct_2$ is a dodecad which may be assumed
to be $O_1 + O_1' + O_{10}$; its complement is also a dodecad. Except for 
these two dodecads, there is no any other dodecad which is a union
of orbits. Thus $N(Rt)^{G_N}$ has a basis:
$$e_i = r_i (i = 1, 2, 3, 4), \,\,
e_5 = \frac{1}{2}(O_1 + O_1' + O_{10}), \,\, e_6 = \frac{1}{2}(O_1'' + O_1''' + O_{10}').$$

\par \noindent
(8) For $Rt = 12A_2$ and the orbit decomposition
$[1, 1, 5, 5, 6, 6]$, 
we label the orbits as $O_1 = \{r_1\}$, $O_1' = \{r_2\}$,
$O_{5} = \{r_3, r_5, \dots, r_{11}\}$,
$O_{5}' = \{r_{4}, r_6, \dots, r_{12}\}$,
$O_{6} = \{r_{13}, r_{15}, \dots, r_{23}\}$,
$O_{6}' = \{r_{14}, r_{16}, \dots, r_{24}\}$,
where $r_{2k-1} + r_{2k}$ ($1 \le k \le 12$)
are the 12 connected components of $Rt$.
Every non-trivial element of the group $N(Rt)/Rt$ is represented by some
$\gamma_H = \sum_{i \in H} \pm (r_{2i-1} + 2 r_{2i})/3$ where
$H$ is an element of the ternary Golay code (so $|H| = 6, 9, 12$) 
which is also the 
Steiner system $St(5, 6, 12)$ [Atlas]. Let $H_i$ with $i = 1$
(resp. $i = 2$) be the unique element of the ternary Golay code with 
$|H_i| = 6$ such that
$\gamma_{H_1} = \frac{1}{3}\sum_{i=2}^6 \pm (r_{2i-1} + 2 r_{2i})
\pm \frac{1}{3}(r_{2j_1-1} + 2 r_{2j_1})$ for some $j_1$
(resp. $\gamma_{H_2} = \frac{1}{3}\sum_{i=7}^{11} \pm (r_{2i-1} + 2 r_{2i})
\pm \frac{1}{3}(r_{2j_2-1} + 2 r_{2j_2})$ for some $j_2$);
such Golay code can also be constructed
from the binary Golay code = Steiner system $St(5, 8, 24)$
where such $H_i$ is the intersection of a fixed dodecad and an octad.
Using the fact that an order-5 element in $A_5$ has
cycle type $(5^2)$ in $M_{12}$ [EDM], as in the case of
$Rt = 24A_1$, we can prove that $N(Rt)^{G_N}$ has a basis:
$$e_i = r_i (i = 1, 2), \,
e_3 = \sum_{k=2}^6 r_{2k-1}, \, e_4 = \sum_{k=7}^{12} r_{2k-1}, \,
e_5 = \frac{1}{3}\sum_{k=1}^6 (r_{2k-1} + 2r_{2k}), \, e_6 = 
\frac{1}{3}\sum_{k=7}^{12} (r_{2k-1} + 2r_{2k}).$$

\par \noindent
(7) For $Rt = 12A_2$ and the orbit decomposition
$[1, 1, 1, 1, 10, 10]$, 
we label the orbits as $O_1 = \{r_1\}$, $O_1' = \{r_2\}$,
$O_{1}'' = \{r_3\}$,
$O_{1}''' = \{r_{4}\}$,
$O_{10} = \{r_{5}, r_{7}, \dots, r_{23}\}$,
$O_{10}' = \{r_{6}, r_{8}, \dots, r_{24}\}$,
where $r_{2k-1} + r_{2k}$ ($1 \le k \le 12$)
are the 12 connected components of $Rt$. As in (8), 
we can prove that $N(Rt)^{G_N}$ has a basis:
$$e_i = r_i \,(1 \le i \le 4), \,\,
e_5 = \sum_{k=3}^{12} r_{2k-1}, \,\, e_6 = 
\frac{1}{3}\sum_{k=1}^{12} (r_{2k-1} + 2r_{2k}).$$

\par \noindent
(9) follows from the direct calculation. Indeed, 
in the case of $M_2$, 
the isomorphism $\varphi_2 : \langle \overline{t}_1^*, \overline{t}_2^* \rangle 
\rightarrow A_{N(Rt)^{G_N}}$
is given by $(\varphi_2(\overline{t}_1^*), \varphi_2(\overline{t}_2^*)) = 
(\overline{e}_1^*, \overline{e}_3^*) \pmatrix 2 & 7 \\ 1 & 0 \endpmatrix$.
In the case of $M_3$, 
the isomorphism $\varphi_3 : \langle \overline{t}_1^*, \overline{t}_2^* \rangle 
\rightarrow A_{N(Rt)^{G_N}}$
is given by $(\varphi_3(\overline{t}_1^*), \varphi_3(\overline{t}_2^*)) = 
(\overline{e}_2^*, \overline{e}_1^* + \overline{e}_4^*) 
\pmatrix 1 & 7 \\ 1 & -4 \endpmatrix$.

\par \noindent
This proves {\bf (2.1)}.

\par \vskip 2pc \noindent
{\bf \S 3. The proof of Theorems A and B}

\par \vskip 1pc \noindent
In this section, we shall prove Theorems A and B.
We prove first the result below which includes Theorem A.

\par \vskip 1pc \noindent
{\bf Theorem 3.1.} 

\par \noindent
(1) There is no faithful group action of the form $A_5 . \mu_3$
(see {\bf (1.0)}) on a $K3$ surface.

\par \noindent
(2) If $G = A_5 . \mu_I$ acts faithfully on a $K3$ surface. Then $G = A_5 : \mu_I$ and
$I = 1$, $2$, or $4$.
(It is proved in [Z2] that $I = 4$ is impossible.)

\par \vskip 1pc \noindent
(2) is a consequence of (1) and {\bf (1.1)}. Indeed, if $I = 6$, 
then the subgroup $H = \alpha^{-1}(\mu_3)$ of $G = A_5 . \mu_6$
is of the form $H = A_5 . \mu_3$ which is impossible by (1).
To prove (1), we need the following result first.

\par \vskip 1pc \noindent
{\bf Lemma 3.2.} Suppose that $G = A_5 . \mu_3$ acts on a $K3$ surface $X$.
Let $\zeta_3 = \text{\rm exp}(2 \pi \sqrt{-1}/3)$.

\par \noindent
(1) We have $G = A_5 \times \mu_3$. Moreover,
a generator $g$ of $\mu_3$ can be chosen so that
$g^*|S_X \otimes {\bold C} = \text{\rm diag}[1, 1, \zeta_3 I_4, \zeta_3^{-1} I_4, \zeta_3 I_5,
\zeta_3^{-1} I_5]$, where the decompositoin here is compatible with
that in {\bf (1.6)} in the sense that $g^* | \chi_4 \oplus \chi_4' = \text{\rm diag}[\zeta_3 I_4, \zeta_3^{-1} I_4]$
and $g^* | \chi_5 \oplus \chi_5' = \text{\rm diag}[\zeta_3 I_5, \zeta_3^{-1} I_5]$.
In particular, $\chi_{\text{\rm top}}(X^g) = -6$.

\par \noindent
(2) We have $S_X^{G} = S_X^{g} = S_X^{A_5} = H^0(X, {\bold Z})^g$. 
This lattice is of rank 2
(whose ${\bold C}$-extension is $\chi_1 \oplus \chi_1'$)
and its discriminant group is 3-elementary.

\par \noindent
(3) We have $S_X^{A_5} = U = U(1)$, or $U(3)$, where
$U(n) = {\bold Z}[u_1, u_2]$ is a rank 2 lattice with
$u_i^2 = 0$ and $u_1 . u_2 = n$.

\par \vskip 1pc \noindent
{\it Proof.} (1) The first part is from {\bf (1.1)}.
For a generator $g$ of $\mu_3$, 
since $o(g) = 3$ and by the form of the decomposition in {\bf (1.6)},
each $\chi_i$ ($i = 4, 5$) is $g$-stable. 
Since the order-3 element $g$ acts on the rank-2 lattice $S_X^{A_5}$
(which is defined over ${\bold Z}$ and
whose ${\bold C}$-extension is $\chi_1 \oplus \chi_1'$),
it has at least one eigenvalue equal to 1
because $G = \langle A_5, g \rangle$ stabilizes an ample line bundle
(the pull back of an ample line bundle on $X/G$).
So $g^* | S_X^{A_5} = \text{\rm id}$.
The commutativity of $g$ with all elements in $A_5$ implies that
$g^* | \chi_i$ is a scalar multiple, by Schur's lemma.

\par \noindent
Thus we can write
$g^* | S_X \otimes {\bold C} = \text{\rm diag} [1, 1, \zeta_3^b I_4,
\zeta_3^c I_4, \zeta_3^d I_5, \zeta_3^e I_5]$, where the ordering 
is the same as in {\bf (1.6)}.
Let $a \in A_5$. Then $(ga)^* | T_X = g^* | T_X$ and the latter
can be diagonalized as $\text{\rm diag}[\zeta_3, \zeta_3^{-1}]$, noting that
$\text{\rm rank} \, T_X = 22 - \text{\rm rank} \, S_X = 2$ 
[Ni1, Theorem 0.1], {\bf (1.0A-B)}.
So $\text{\rm Tr}(g a)^* | T_X = -1$.
As in the proof of {\bf (1.8)}, the topological Lefschetz fixed point
formula implies that
$\chi(X^{ga}) = 2 + \text{\rm Tr}  (ga)^*|T_X + \text{\rm Tr}(ga)^*|S_X = 1 + \text{\rm Tr}(ga)^*|S_X =
3 + \zeta^b \text{\rm Tr}(a^* | \chi_4) + \zeta^c \text{\rm Tr}(a^* | \chi_4') +  
\zeta^d \text{\rm Tr}(a^* | \chi_5) + \zeta^e \text{\rm Tr}(a^* | \chi_5')$.
So for $a = \text{\rm id}$, $2A$, $3A$, $5A$ with $nA$ denoting an element of order
$n$ in $A_5$, we have:
$$\gather 
\chi_{\text{\rm top}}(X^g) = 3 + 4(\zeta_3^b + \zeta_3^c) + 
5(\zeta_3^d + \zeta_3^e), \\
\chi_{\text{\rm top}}(X^{g2A}) = 3 + \zeta_3^d + \zeta_3^e, \\
\chi_{\text{\rm top}}(X^{g3A}) = 3 + \zeta_3^b + \zeta_3^c  
-\zeta_3^d - \zeta_3^e, \\
\chi_{\text{\rm top}}(X^{g5A}) = 3 - \zeta_3^b - \zeta_3^c.
\endgather$$

\par \noindent
The fact $\chi(X^{g5A}) = 4$ in {\bf (1.4)} implies that
$(\zeta_3^b, \zeta_3^c) = (\zeta_3, \zeta_3^{-1})$ after switching
$\chi_4$ with $\chi_4'$ if necessary.
Since $\chi(X^{g 3A}) = 0$ is in ${\bold R}$ (in ${\bold Z}$, indeed),
we may assume that $(\zeta_3^d, \zeta_3^e) = (\zeta_3, \zeta_3^{-1})$,
or $(1, 1)$. If the former case occurs then the lemma is true.

\par \noindent
Suppose that the latter case occurs. Then $\chi_{\text{\rm top}}(X^g) = 9$, whence
$n_g = 2$ and $|X^g_{\text{\rm isol}}| = m_g = n_g + 3 = 5$ by {\bf (1.4)}.
Since $g$ commutes with every element in $A_5$, our $A_5$
acts on the 5-point set $X^g_{\text{\rm isol}}$. By {\bf (1.7)},
$A_5$ either fixes a point $P_1$ of the set (and hence $A_5 < SL(T_{X, P_1})$,
contradicting {\bf (1.0C)}), or acts transitively as a subgroup ($= [S_5, S_5]$)
of $S_5$, on the set
with an order-12 stabilizer (of a point $P_1$) subgroup $A_4 < A_5$,
so $A_4 < SL(T_{X, P_1})$, contradicting
{\bf (1.0C)}.
This proves the assertion (1).

\par \noindent
(2) The first part follows from (1), that $g^* | T_X \otimes {\bold C} = 
\text{\rm diag}[\zeta_3, \zeta_3^{-1}]$ w.r.t. to a suitable basis by [Ni1, Theorem 0.1]
and that all lattices in (2)
are primitive (of the same rank as they turn out to be)
in $L := H^2(X, {\bold Z})$. We still have to show that
the discriminant group $A_{L^g} = \text{\rm Hom}(L^g, {\bold Z})/L^g$ 
of $L^g$ is 3-elementary.
Let $L_g = (L^g)^{\perp}$ be the orthogonal of $L^g$ in $L$.
Then $g^* | L_g$ has only $\zeta_3^{\pm}$ as eigenvalues. Now arguing
as in [OZ2, Lemma (1.3)]
(for the finite index sublattice $L^g \oplus L_g$ of $L$,
instead of $S_X \oplus T_X$), we can show that $A_{L^g}$ is
3-elementary.

\par \noindent
(3) follows from (2). See [CS, Table 15.2a].

\par \vskip 1pc \noindent
The fixed locus $X^g$ can be determined:

\par \vskip 1pc \noindent
{\bf Lemma 3.3.} (1) With the assumption and notation in {\bf (3.1)} 
and {\bf (3.2)},
either $X^g = C \coprod R$ is a disjoint union of a genus-5 curve $C$
and a curve $R$ ($\cong {\bold P}^1$) (so $C^2 = 8$, and $S_X^g = U
\supset {\bold Z}[C, R]$), or
$X^g$ equals a sinlge genus-4 curve $C$ (so $C^2 = 6$).

\par \noindent
(2) In the former case, $\Phi_{|C|} : X \rightarrow {\bold P}^5$
is a degree-2 morphism onto either the Veronese-embedded 
${\bold P}^2$ in ${\bold P}^5$
or the normal cone $\overline{\Sigma}_4$ over a rational normal twisted quartic in ${\bold P}^4$.

\par \vskip 1pc \noindent
{\it Proof.} Since $\chi(X^g) = -6$ by {\bf (3.2)}, we have 
$n_g = -3$ and $m_g = 0$ in notation of {\bf (1.4)}.
$n(g) < 0$ infers that $X^g$ is a disjoint union of a smooth curve $C$ of genus $\ge 2$ and $t$ of ${\bold P}^1$'s with $-6 = 2 - 2g(C) + 2t$ (see
{\bf (1.2)}). The fact that
$\text{\rm rank} \, S_X^{g} = 2$ in {\bf (3.2)} implies that either $t = 0$ (so $g(C) = 4$),
or $t = 1$ (so $g(C) = 5$) 
so that the two curves in $X^g$ give rise to two linearly
independent classes of $S_X^g$.

\par \noindent
If $S_X^g = U(3)$, then $C^2 = 0$ (mod 3) because $C$ is in $S_X^g$,
whence $C^2 = 6$. This proves the first assertion of the lemma, by virtue of {\bf (3.2)}.

\par \noindent
Consider the caser $X^g = C \coprod R$.
By [SD, Theorem 3.1], $|C|$ is base point free and we have
a morphism $\varphi := \Phi_{|C|} : X \rightarrow {\bold P}^5$.
Now $8 = C^2 = \text{\rm deg}( \varphi) . \text{\rm deg}( \text{\rm Im} \, \varphi)$,
where $\text{\rm deg}( \text{\rm Im} \, \varphi) \ge 5-1$. Thus either $\varphi$ is an
embedding modulo the curves in $C^{\perp}$, or
$\varphi$ is a degree-2 map as described in {\bf (3.3)} [SD, Theorem 5.2, 
Propositions 5.6 and 5.7].

\par \noindent
Write $S_X^g = {\bold Z}[u_1, u_2]$ with $u_i^2 = 0$ and $u_1 . u_2 = 1$.
Express $C = a_1 u_1 + a_2 u_2$. Then $8 = 2a_1 a_2$ and we may assume that
$(a_1, a_2) = (2, 2)$ or 
$(4, 1)$ (after replacing $u_i$ by $-u_i$ or switching $u_1$
with $u_2$ if necessary). So $C . u_i > 0$ and hence 
the Riemann-Roch theorem implies that $\dim |u_1| \ge 1$.
Write $|u_1| = |M| + F$ with $|M|$ the movable part.
Then $0 < C . M \le C . u_1 = a_2 \le 2$. If
$\varphi$ is birational then 
$\varphi(M)$ is a plane conic or a line, whence
$M \cong {\bold P}^1$, $M^2 = -2$ and $|M|$ is not movable,
a contradiction. This proves the lemma.

\par \vskip 1pc \noindent
We now prove {\bf (3.1)} (1). 
Consider the case in {\bf (3.3)}, where $X^g = C \coprod R$ 
and $\varphi = \Phi_{|C|} : X \rightarrow {\bold P}^5$ is 
a degree-2 morphism onto the Veronese-embedded ${\bold P}^2$
in ${\bold P}^5$. Since $C$ (and hence $|C|$) is $G$-stable, there is an induced
action of $G$ on ${\bold P}^5$ (and hence also an action
of $G$ on the image $\varphi(X) = {\bold P}^2$)
so that the map $\varphi$ is $G$-equivariant. 
The $G = A_5 \times \mu_3$ action on the image is also faithful because 
$A_5$ is simple and $\text{\rm deg}( \varphi) = 2$ is coprime to $3$
($= |\mu_3|$).
The action of $A_5$ on the image is via $A_5 \subset SL_3({\bold C}) \subset
PGL_3({\bold C})$ and is given in Burnside [Bu, \S 232, or \S 266] {\bf (1.8).}
In particular, the commutativity of $g$ with the two generators (order 5 and 2)
of $A_5$ in [Bu, \S 266] shows that $g$ is a scalar and acts trivially
on the image $\varphi(X) = {\bold P}^2$, a contradiction.

\par \noindent
Consider the case in {\bf (3.3)}, where $X^g = C \coprod R$ in {\bf (3.3)}
and $\varphi = \Phi_{|C|} : X \rightarrow {\bold P}^5$ is 
a degree-2 morphism onto the cone $\overline{\Sigma}_4$.
Note that the minimal resolution $\Sigma_4$ of $\overline{\Sigma}_4$
is the Hirzebruch surface of degree 4. As in the previous case,
there is a faithful action of $G$ on $\overline{\Sigma}_4$ such 
that $\varphi$ is $G$-equivariant. Note that the image $\varphi(C)$
is a hyperplane section away from the singularity and with
$\varphi(C)^2 = 4$. Let $\ell$ be a generating line of the cone
$\overline{\Sigma}_4$. Then $\varphi(C) \sim 4 \ell$ as Weil divisors.
This gives rise to a ${\bold Z}/(4)$-cover
$\pi : Y = Spec \oplus_{i=0}^{3} {\Cal O}_{\overline{\Sigma}_4}(-i \ell)
\rightarrow \overline{\Sigma}_4$ which is
(totally) ramified exactly over $\varphi(C)$. One sees that $Y \cong {\bold P}^2$
and $\pi^* \varphi(C) = 4 L$ with $L$ a line.
Clearly, $A_5$ ($< G$) stabilizes the divisorial sheaves ${\Cal O}(-i \ell)$
and fixes the defining equation of $\varphi(C)$, so there is an
induced faithful $A_5$-action on $Y = {\bold P}^2$ so that
$\pi$ is $A_5$-equivariant (see {\bf (1.7)}). Now $L$ is stabilized by $A_5$
(because so is $\varphi(C)$).
So the defining equation $F_1 = 0$ of $L$ is semi $A_5$-invariant
(and hence $A_5$-invariant because of the simplicity of the group $A_5$).
But every $A_5$-invariant form is of even degree by [Bu, \S 266],
noting also that the action of $A_5$ on ${\bold P}^2$ is via
$A_5 \subset SL_3({\bold C}) \rightarrow PGL_3({\bold C})$ by {\bf (1.8)}.
We reach a contradiction.

\par \noindent
Consider the case $X^g = C$ in {\bf (3.3)}. 
Let $f : X \rightarrow Y = X/\langle g \rangle$
be the quotient map. There is an induced faithful action
$A_5$ on $Y$ so that $f$ is $A_5$-equivariant.
Then by the ramification divisor formula,
$0 \sim K_X = f^*(K_Y) + 2C$. Pushing down by $f_*$,
one obtains $0 \sim 3K_Y + 2 B$ with $B = f_*C = f(C) \cong C$
and $f^*B = 3C$, so $B^2 = 3 C^2 = 18$.
Solving, one obtains $B = (-3/2)K_Y$ and $K_Y^2 = 8$.
Thus the smooth ruled surface $Y$ equals 
a Hirzebruch surface $\Sigma_d$ of degree $d$.
The irreducibility of $B$ (being a {\bf Z}-divisor)
implies that $d = 0, 2$ [Ha, Ch V, Cor. 2.18].

\par \noindent
Suppose that $d = 2$. Then the (-2)-curve $M$
on $Y$ is disjoint from $B = (-3/2)K_Y$
and hence $f^*M = \coprod_{i=1}^3 M_i$ is a disjoint union
of three (-2)-curves not intersecting $C$.
Since $M$ is clearly $A_5$-stable, the set $\coprod M_i$ is also
$A_5$-stable, whence each $M_i$ is $A_5$-stable {\bf (1.7)}.
An order-5 element $5A$ in $A_5$ acts on each $M_i$  faithfully
by {\bf (1.2)} and has exactly two fixed points by {\bf (1.8)}.
But according to {\bf (1.2)}, 
$4 = |X^{5A}| \ge \sum_{i=1}^3 |M_i^{5A}| = 6$, a contradiction.

\par \noindent
Thus $d = 0$. Clearly, the simple group $A_5$ stabilizes each ruling
and there is an induced action $\rho_i : A_5 \times {\bold P}^1 
\rightarrow {\bold P}^1$ with $i = 1, 2$
for the $i$-th ${\bold P}^1$ in $Y = \Sigma_0 = {\bold P}^1 \times {\bold P}^1$
so that $\rho_1 \times \rho_2$ is the given $A_5$ action on $Y$.
Changing coordinates suitably we may assume that the $A_5$ action on $Y$
commutes with the involution $\iota$ of $Y$ switching the  
two ${\bold P}^1$'s in $Y$, i.e., $\rho_1 = \rho_2$ as 
actions of $A_5$ on the same ${\bold P}^1$ {\bf (1.8)}.

\par \noindent
Let $j : Y \rightarrow Z = Y/\langle \iota \rangle = {\bold P}^2$
be the quotient map. Then there is an induced faithful action of $A_5$
on ${\bold P}^2$ such that $j$ is $A_5$-equivariant. Now 
$\iota(B)$ is an irreducible curve with $(\iota(B))^2 = B^2/2 = 9$.
It is a cubic curve and $A_5$-stable
because so is $B = f(C)$.
The action of $A_5$ on $Y/\langle \iota \rangle = {\bold P}^2$
is via $SL_3({\bold C}) \rightarrow PGL_3({\bold C})$ {\bf (1.8)}.
The defining equation $F_3$ of $\iota(B)$ is then a cubic form
and semi $A_5$-invariant (and hence $A_5$-invariant by the
simplicity of the group $A_5$). However, Burnside [Bu, \S 266]
shows that every $A_5$-invariant form is of even degree, a contradiction.
This completes the proof of {\bf (3.1)} (1) and also of {\bf (3.1)}.

\par \vskip 1pc \noindent
We now prove Theorem B.
Suppose that $G = A_5 . C_n$ acts faithfully on a $K3$ surface $X$.
By {\bf (1.0A)}, $A_5 \le G_N$. So $G_N = A_5$, $S_5$, $A_6$ or 
$M_{20} = C_2^{\oplus 4} : A_5$
by [Xi, the list]. In notation of {\bf (1.0)}, for some $m \,|\, n$, we have
$G_N = \text{\rm Ker}(\alpha) = A_5 . C_m$ 
and $G/G_N = \mu_I$, where $n = mI$. 
By the same proof of {\bf (1.1)}, we have $I = 1, 2, 3, 4$, or $6$.
Let $h \in G$ such that the coset of $h$ is a generator
of $G/A_5 = C_n$. Then $h^* \omega_X = \eta_I \omega_X$ for some
primitive $I$-th root $\eta_I$ of $1$.
Note that $n \,|\, \text{\rm ord}(h)$ and $h^I \in G_N$, whence $\text{\rm ord}(h^I) \le 8$
by {\bf (1.2)}. Thus $\text{\rm ord}(h) = I \, \text{\rm ord}(h^I)$ and $m \,|\, \text{\rm ord}(h^I)$.
In particular, $|G_N| = m |A_5| \le 8 |A_5|$.
Hence $G_N \ne M_{20}$.

\par \noindent
If $G_N = A_5 . C_m = A_6$, then $m = 6$ and $A_6$ 
includes $\langle h^I \rangle \ge C_6$, which is impossible.
If $G_N = A_5$, then $C_n = \mu_I$ in notation of {\bf (1.0)},
and Theorem B follows from {\bf (3.1)}.

\par \noindent
Consider the case $G_N = S_5$. Then $m = 2$ and $n = 2I$.
Moreover, $G_N = \langle A_5, h^I \rangle$. So $h^I \in S_5 - A_5$.
Since $S_5 \rightarrow \text{\rm Aut}(S_5)$ ($x \mapsto c_x$) is an isomorphism,
we have $c_h = c_s$ for some $s \in S_5$. Set $g = h s^{-1}$.
Then $g$ commutes with every element in $S_5$ and also $\alpha(g) = \alpha(h)$
is a generator of $\text{\rm Im} (\alpha) = \mu_I$. Now
$g^I \in \text{\rm Ker}(\alpha) = G_N = S_5$ is in the centre of $S_5$ (which is $(1)$).
So $\text{\rm ord}(g) = I$ and $G = G_N \times \langle g \rangle = S_5 \times \mu_I
> A_5 \times \mu_I$. By {\bf (3.1)}, we have $I = 1, 2$ or $4$.
The $S_5 \times \mu_I$ should have an element $h$ such that $h^I \in S_5 - A_5$
(i.e., $h^I$ is not an even permutation). Thus, $I \ne 2$, or $4$.
Therefore, $I = 1$ and $G = G_N = S_5$. This completes the proof of Theorem B.

\par \vskip 2pc \noindent
{\bf References}

\par \noindent
[AS2] M. F. Atiyah and G. B. Segal,
The index of elliptic operators. II. 
Ann. of Math. 87 (1968), 531--545.

\par \noindent
[AS3] M. F. Atiyah and I. M. Singer,
The index of elliptic operators. III. 
Ann. of Math. 87 (1968) 546--604.

\par \noindent
[Atlas] J. H. Conway, R. T. Curtis, S. P. Norton, R. A. Parker and R. A. Wilson, 
Atlas of finite groups. Oxford University Press. Reprinted 2003 (with corrections). 

\par \noindent
[Bu] W. Burnside, Theory of groups of finite order. 
Dover Publications, Inc., New York, 1955.

\par \noindent
[CS] J. H. Conway and N. J. A. Sloane, Sphere packings, lattices and groups. 3rd ed. 
Grundlehren der Mathematischen Wissenschaften, 290. Springer-Verlag, New York, 1999.

\par \noindent
[EDM] Encyclopedic dictionary of mathematics. Vol. I--IV. Translated from the Japanese. 
2nd ed. Edited by Kiyosi Itô. MIT Press, Cambridge, MA, 1987.

\par \noindent
[Ha] R. Hartshorne, Algebraic geometry. 
Graduate Texts in Mathematics, No. 52. 
Springer - Verlag, New York-Heidelberg, 1977.

\par \noindent
[IOZ] A. Ivanov, K. Oguiso and D. -Q. Zhang, The monster and $K3$ surfaces, in preparation.

\par \noindent
[KOZ1] J. Keum, K. Oguiso and D. -Q. Zhang, The alternating group of degree 6
in geometry of the Leech lattice and $K3$ surfaces, 
math.AG/0311462, Proc. London Math. Soc. to appear.

\par \noindent
[KOZ2] J. Keum, K. Oguiso and D. -Q. Zhang, Extensions of 
the alternating group of degree 6 in geometry of $K3$ surfaces,
math.AG/0408105.

\par \noindent
[Ko1] S. Kondo, Niemeier lattices, Mathieu groups, and finite groups of symplectic 
automorphisms of $K3$ surfaces. 
Duke Math. J. 92 (1998), 593--598.

\par \noindent
[Ko2] S. Kondo, The maximum order of finite groups of automorphisms of $K3$ surfaces.
Amer. J. Math. 121 (1999), 1245--1252.

\par \noindent
[MO] N. Machida and K. Oguiso,
On $K3$ surfaces admitting finite non-symplectic group actions.  
J. Math. Sci. Univ. Tokyo 5 (1998), 273--297.

\par \noindent
[Mu1] S. Mukai, Finite groups of automorphisms of $K3$ surfaces and the Mathieu group. 
Invent. Math. 94 (1988), 183--221.

\par \noindent
[Mu2] Lattice-theoretic construction of symplectic actions on $K3$ surfaces,
Duke Math. J. 92 (1998), 599--603. As the Appendix to [Ko1].

\par \noindent
[Ni1] V. V. Nikulin, Finite automorphism groups of Kahler $K3$ surfaces,
Trans. Moscow Math. Soc. 38 (1980), 71--135.

\par \noindent
[Ni2] V. V. Nikulin, Factor groups of groups of automorphisms of hyperbolic forms
with respect to subgroups generated by 2-reflections. Algebrogeometric applications.
J. Soviet Math. 22 (1983), 1401--1475.

\par \noindent
[Og] K. Oguiso, A characterization of the Fermat quartic $K3$ surface by
means of finite symmetries, math.AG/0308062. Compositio Math. to appear.

\par \noindent
[OZ1] K. Oguiso and D. -Q. Zhang, 
On the most algebraic $K3$ surfaces and the most extremal log Enriques surfaces. Amer. J. Math. 118 (1996), 1277--1297.

\par \noindent
[OZ2] K. Oguiso and D. -Q. Zhang,
On Vorontsov's theorem on $K3$ surfaces with non-symplectic group actions. 
Proc. Amer. Math. Soc. 128 (2000), 1571--1580.

\par \noindent
[OZ3] K. Oguiso and D. -Q. Zhang, 
The simple group of order 168 and $K3$ surfaces.
Complex geometry (Gottingen, 2000), Collection of papers
dedicated to Hans Grauert, 165--184, 
Springer, Berlin, 2002. 

\par \noindent
[SD] B. Saint-Donat,
Projective models of $K-3$ surfaces. 
Amer. J. Math. 96 (1974), 602--639.

\par \noindent
[Su] M. Suzuki, Group theory. I. 
Translated from the Japanese by the author. 
Grundlehren der Mathematischen Wissenschaften 247. 
Springer-Verlag, Berlin-New York, 1982.

\par \noindent
[Xi] G. Xiao, Galois covers between $K3$ surfaces.
Ann. Inst. Fourier (Grenoble) 46 (1996), 73--88.

\par \noindent
[Z1] D. -Q. Zhang, Quotients of $K3$ surfaces modulo involutions. 
Japan. J. Math. (N.S.) 24 (1998), 335--366.

\par \noindent
[Z2] D. -Q. Zhang, The alternating groups and K3 surfaces, preprint 2004.

\par \vskip 2pc \noindent
D. -Q. Zhang
\par \noindent
Department of Mathematics
\par \noindent
National University of Singapore
\par \noindent
Singapore
\par \noindent
E-mail : MATZDQ$\@$MATH.NUS.EDU.SG

\enddocument